\documentclass{amsart}

\usepackage{amsmath,amsthm,amssymb}
\usepackage{hyperref}
\usepackage{pstricks,pst-node,multido}

\begin{document}

\title{A note on 2-distant noncrossing partitions and weighted Motzkin paths}

\author{Ira M. Gessel}
\author{Jang Soo Kim}
\thanks{The first author was supported by NSA Grant H98230-10-1-0196. The second
  author was supported by the grant ANR08-JCJC-0011.} 
\email[Ira M. Gessel]{gessel@brandeis.edu}
\email[Jang Soo Kim]{kimjs@math.umn.edu}

\begin{abstract}
  We prove a conjecture of Drake and Kim: the number of $2$-distant noncrossing
  partitions of $\{1,2,\ldots,n\}$ is equal to the sum of weights of Motzkin
  paths of length $n$, where the weight of a Motzkin path is a product of
  certain fractions involving Fibonacci numbers. We provide two proofs of their
  conjecture: one uses continued fractions and the other is combinatorial.
\end{abstract}

\newtheorem{thm}{Theorem}[section]
\newtheorem{lem}[thm]{Lemma}
\newtheorem{prop}[thm]{Proposition}
\newtheorem{cor}[thm]{Corollary}
\theoremstyle{definition}
\newtheorem{example}{Example}
\newtheorem{defn}{Definition}
\newtheorem{conj}{Conjecture}
\newtheorem{question}{Question}
\theoremstyle{remark}
\newtheorem{remark}{Remark}
\newtheorem*{note}{Note}
\def\sch.{Schr\"oder}
\newcommand{\ds}{\displaystyle}
\renewcommand{\aa}{\alpha}
\newcommand{\bb}{\beta}
\newcommand{\5} {\sqrt{(1-x)(1-5x)}}
\newcommand{\NC} {\operatorname{NC}}
\newcommand{\scheven} {\operatorname{SCH}_{\rm even}}
\newcommand{\dyck} {\operatorname{Dyck}}
\newcommand{\mot} {\operatorname{Mot}}

\psset{unit=1cm, dotsize=5pt}
\psset{subgriddiv=1,gridcolor=lightgray,gridlabels=0pt,gridwidth=0.4pt, linewidth=1.5pt}

\newcount\sx \newcount\sy \newcount\ex \newcount\ey
\def\CHS#1(#2,#3){ 
\sx=#2 \sy=#3 \ex=#2 \ey=#3
\advance\ex by1 \advance\ey by0
\psline(\sx,\sy)(\ex,\ey)
\rput(\number\sx.5,\number\sy.3){$#1$}
\psdot(\number\sx,\number\sy) \psdot(\number\ex,\number\ey)
}
\def\CDS#1(#2,#3){ 
\sx=#2 \sy=#3 \ex=#2 \ey=#3
\advance\ex by1 \advance\ey by-1
\psline(\sx,\sy)(\ex,\ey)
\rput(\number\sx.8,\number\ey.8){$#1$}
\psdot(\number\sx,\number\sy) \psdot(\number\ex,\number\ey)
}
\def\HS(#1,#2){ 
\sx=#1 \sy=#2 \ex=#1 \ey=#2
\advance\ex by1 \advance\ey by0
\psline(\sx,\sy)(\ex,\ey)
\psdot(\number\sx,\number\sy) \psdot(\number\ex,\number\ey)
}
\def\US(#1,#2){ 
\sx=#1 \sy=#2 \ex=#1 \ey=#2
\advance\ex by1 \advance\ey by1
\psline(\sx,\sy)(\ex,\ey)
\psdot(\number\sx,\number\sy) \psdot(\number\ex,\number\ey)
}
\def\DS(#1,#2){ 
\sx=#1 \sy=#2 \ex=#1 \ey=#2
\advance\ex by1 \advance\ey by-1
\psline(\sx,\sy)(\ex,\ey)
\psdot(\number\sx,\number\sy) \psdot(\number\ex,\number\ey)
}
\def\DHS(#1,#2){ 
\sx=#1 \sy=#2 \ex=#1 \ey=#2
\advance\ex by2 \advance\ey by0
\psline(\sx,\sy)(\ex,\ey)
\psdot(\number\sx,\number\sy) \psdot(\number\ex,\number\ey)
}


\maketitle

\section{Introduction}

A \emph{Motzkin path} of length $n$ is a lattice path from $(0,0)$ to $(n,0)$
consisting of up steps $U=(1,1)$, down steps $D=(1,-1)$ and horizontal steps
$H=(1,0)$ that never goes below the $x$-axis. The \emph{height} of a step in a
Motzkin path is the $y$ coordinate of the ending point.

\begin{figure}
  \centering
\begin{pspicture}(0,0)(8,2) 
\psgrid(0,0)(8,2) 
\US(0,0) \CDS{b_0}(1,1) \CHS{\lambda_0}(2,0) \US(3,0) \US(4,1) 
\CDS{b_1}(5,2) \CHS{\lambda_1}(6,1) \CDS{b_0}(7,1)
\end{pspicture}
\caption{A Motzkin path and the weights of its steps with respect to
  $(b,\lambda)$.}
  \label{fig:mot}
\end{figure}
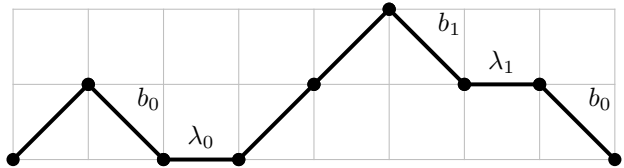

Given two sequences $b=(b_0,b_1,\ldots)$ and
$\lambda=(\lambda_0,\lambda_1,\ldots)$, the \emph{weight} of a Motzkin path with
respect to $(b,\lambda)$ is the product of $b_i$ and $\lambda_i$ for each
horizontal step and down step of height $i$ respectively, see
Figure~\ref{fig:mot}. Let $\mot_n(b,\lambda)$ denote the sum of weights of
Motzkin paths of length $n$ with respect to $(b,\lambda)$. This sum is closely
related to orthogonal polynomials; see \cite{Viennot1983, Viennot1985}.

Drake and Kim \cite{Drake2008} defined the set $\NC_k(n)$ of $k$-distant
noncrossing partitions of $[n]=\{1,2,\ldots,n\}$. For $k\geq0$, a
\emph{$k$-distant noncrossing partition} is a set partition of $[n]$ without two
arcs $(a,c)$ and $(b,d)$ satisfying $a<b\leq c<d$ and $c-b\geq k$, where an
\emph{arc} is a pair $(i,j)$ of integers contained in the same block which does
not contain any integer between them. For example,
$\pi=\{\{1,5,7\},\{2,3,6\},\{4\}\}$ is a $3$-distant noncrossing partition but
not a $2$-distant noncrossing partition because $\pi$ has two arcs $(1,5)$ and
$(3,6)$ with $5-3\geq2$. Note that the $1$-distant noncrossing partitions are
the ordinary noncrossing partitions, which implies that $\#\NC_1(n)$ is equal to
the Catalan number $\frac{1}{n+1}\binom{2n}{n}$. It is not difficult to see that
$\NC_0(n)$ is in bijection with the set of Motzkin paths of length $n$. In the
same paper, they proved that
\begin{equation}
  \label{eq:gf}
\sum_{n\geq0} \#\NC_2(n) x^n = \frac32 -\frac12 \sqrt{\frac{1-5x}{1-x}}.
\end{equation}
The number $\#\NC_2(n)$ also counts many combinatorial objects: \sch. paths with
no peaks at even levels, etc; see \cite{Kim2009b,Mansour2007, Yan2009}.

There are simple expressions of $\#\NC_k(n)$ using Motzkin paths for $k=0,1,3$:
\begin{align*}
\#\NC_0(n) & =\mot_n ((1,1,\ldots),(1,1,\ldots)), \\
\#\NC_1(n)&=\mot_n ((1,2,2,\ldots),(1,1,\ldots)),\\
\#\NC_3(n)&=\mot_n ((1,2,3,3,\ldots),(1,2,2,\ldots)),
\end{align*}
where the second equation is well known and the third one was first conjectured
by Drake and Kim \cite{Drake2008} and proved by Kim \cite{Kim_front}.
The main purpose of this paper is to prove the following theorem which was also
conjectured by Drake and Kim \cite{Drake2008}.

\begin{thm}\label{thm:main}
  Let $b=(b_0,b_1,\ldots)$ and $\lambda=(\lambda_0,\lambda_1,\ldots)$ be the
  sequences with $b_0=\lambda_0=1$ and for $n\geq1$,
\begin{equation}
  \label{eq:3}
 b_n = 3-\frac{1}{F_{2n-1}F_{2n-3}} \quad \text{and} \quad
  \lambda_n = 1+\frac1{F_{2n-1}^2},
\end{equation}
where $F_m$ is the Fibonacci number defined by $F_0=0, F_1=1$, and
$F_{m}=F_{m-1}+F_{m-2}$ for all $m$ (so $F_{-1}=1$). Then we have
\[\#\NC_2(n)=\mot_n (b,\lambda).\]
\end{thm}

Theorem~\ref{thm:main} is very interesting because it is not even obvious that
$\mot_n (b,\lambda)$ is an integer.  In this paper, we give two proofs of
Theorem~\ref{thm:main}: one uses continued fractions and the other is
combinatorial.

\section{Continued fractions}

Let $\aa=(\aa_0,\aa_1,\aa_2,\ldots)$, $\bb=(\bb_0,\bb_1,\bb_2,\ldots)$, and $c =
(c_0,c_1,c_2, \ldots)$ be sequences of numbers. 

Let $J(x; \aa_0,\bb_0; \aa_1,\bb_1; \aa_2,\bb_2;\dots)= J(x; \aa,\bb)$ denote the
\emph{J-fraction}
\[
\cfrac{1}{ 1-\aa_0 x - 
\cfrac{\bb_0 x^2}{ 1-\aa_1 x -
\cfrac{\bb_1 x^2}{ 1-\aa_2 x - \dots}}}
\]
and let $S(x; c_0, c_1, c_2,\dots)=S(x; c)$ denote the \emph{S-fraction}
\[\cfrac{1}{1-
\cfrac{\mathstrut c_0x}{1-
\cfrac{\mathstrut c_1x}{1-\dots
}}}.
\]

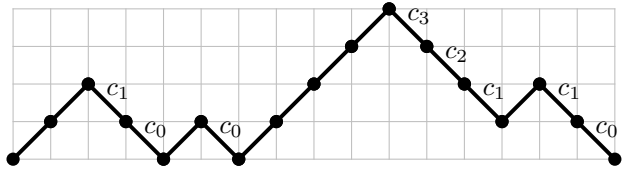
\begin{figure}
  \centering
\psset{unit=0.5cm}
\begin{pspicture}(0,0)(16,4) 
\psgrid(0,0)(16,4) 
\US(0,0) \US(1,1) \CDS{c_1}(2,2) \CDS{c_0}(3,1) \US(4,0) \CDS{c_0}(5,1)
\US(6,0) \US(7,1) \US(8,2) \US(9,3) \CDS{c_3}(10,4) \CDS{c_2}(11,3)
\CDS{c_1}(12,2) \US(13,1) \CDS{c_1}(14,2) \CDS{c_0}(15,1)
\end{pspicture}
 \caption{A Dyck path and the weights of its steps with respect to $c$.}
  \label{fig:dyck}
\end{figure}

A \emph{Dyck path} of length $2n$ is a lattice path from $(0,0)$ to $(2n,0)$
consisting of up steps $U=(1,1)$ and down steps $D=(1,-1)$ that never goes below
the $x$-axis. The \emph{height} of a step in a Dyck path is the $y$ coordinate
of the ending point. The \emph{weight} of a Dyck path with respect to $c$ is the
product of $c_i$ for each down step of height $i$, see Figure~\ref{fig:dyck}.
Let $\dyck_n(c)$ denote the sum of weights of Dyck paths of length $2n$ with
respect to $c$.

It is well known that
\[ \sum_{n\geq0} \mot_n(\aa,\bb) x^n = J(x; \aa,\bb)  \text{\quad and \quad}
\sum_{n\geq0} \dyck_n(c) x^n = S(x;c).\]
The following proposition is easy to see. 

\begin{prop}\label{thm:1}
  If $\aa_n = c_{2n-1}+c_{2n}$ and $\bb_n = c_{2n}c_{2n+1}$ for all $n\geq0$,
  with $c_{-1}=0$, then $S(x; c) = J(x; \aa,\bb)$.
\end{prop}

One can prove Proposition~\ref{thm:1} by the following observation: a Motzkin
path may be obtained from a Dyck path by taking steps two at a time and changing
$UU$, $UD$, $DU$ and $DD$, respectively, to $U$, $H$, $H$ and $D$. For example,
the Motzkin path in Figure~\ref{fig:mot} is obtained from the Dyck path in
Figure~\ref{fig:dyck} in this way.

Let $d=(d_0,d_1,d_2,\dots)$ be the sequence with $d_0=1$ and for $n\geq1$, 
\begin{equation}
  \label{eq:7}
d_{2n-1} = \frac{F_{2n-1}}{F_{2n-3}}, \qquad d_{2n} = \frac{1}{d_{2n-1}}.  
\end{equation}

Recall the two sequences $b=(b_0,b_1,\ldots)$ and
$\lambda=(\lambda_0,\lambda_1,\ldots)$ defined in \eqref{eq:3}. 

\begin{lem}\label{thm:2}
We have the following.
  \begin{enumerate}
  \item $b_n = d_{2n-1} + d_{2n}$ for all $n\geq0$, where $d_{-1}=0$. 
  \item $\lambda_n = d_{2n} d_{2n+1}$ for all $n\geq0$. 
  \item $1/d_{2n-1} + d_{2n+1} = 3$ for all $n\geq1$. 
 \end{enumerate}
\end{lem}
\begin{proof}
We will use two cases of the well-known Catalan identity for Fibonacci numbers, $F_m^2 - F_{m+i}F_{m-i}=(-1)^{m-i}F_i^2$.

  (1) This is true for $n=0$. For $n\geq1$ we have
 \begin{align*}
  d_{2n-1} + d_{2n} &=
  \frac{F_{2n-1}}{F_{2n-3}} + \frac{F_{2n-3}}{F_{2n-1}}=
  \frac{F_{2n-1}^2+F_{2n-3}^2}{F_{2n-1}F_{2n-3}}=
  \frac{2F_{2n-1}F_{2n-3}+(F_{2n-1}-F_{2n-3})^2}{F_{2n-1}F_{2n-3}}
  \\
  &= 
  2+ \frac{F_{2n-2}^2}{F_{2n-1}F_{2n-3}} = 
  3+ \frac{F_{2n-2}^2-F_{2n-1}F_{2n-3}}{F_{2n-1}F_{2n-3}} =
  3-
  \frac{1}{F_{2n-1}F_{2n-3}}=b_n.
  \end{align*}

  (2) This is true for $n=0$. For $n\ge 1$ we have
  \begin{equation*}
  d_{2n}d_{2n+1} =
  \frac{F_{2n-3}}{F_{2n-1}} \frac{F_{2n+1}}{F_{2n-1}} =
  \frac{F_{2n-1}^2 + (F_{2n-3}F_{2n+1} -F_{2n-1}^2)}{F_{2n-1}^2}=
   1+\frac1 {F_{2n-1}^2}=\lambda_n.
   \end{equation*}

(3) We have
\begin{align*}
\frac{1}{d_{2n-1}} + d_{2n+1} &= \frac{F_{2n-3}}{F_{2n-1}} +
\frac{F_{2n+1}}{F_{2n-1}} = \frac{(F_{2n-1} - F_{2n-2})+(F_{2n}+F_{2n-1})}{F_{2n-1}} \\
&=
2 + \frac{F_{2n} - F_{2n-2}}{F_{2n-1}}=3.
\end{align*}

\end{proof}

By Proposition~\ref{thm:1} and Lemma~\ref{thm:2}, we obtain the following.
\begin{cor}\label{thm:4}
For the sequences $b$, $\lambda$ and $d$ defined in \eqref{eq:3} and
\eqref{eq:7}, we have
\[\dyck_n(d) = \mot_n(b,\lambda).\]
\end{cor}

Now we can prove the following $S$-fraction formula for the generating function
\eqref{eq:gf} for $\#\NC_2(n)$.

\begin{thm}
\label{t-1}
We have
\begin{align*}
\frac32 - \frac12\sqrt{\frac{1-5x}{1-x}}
&=S(x; 1,1,1,2,\tfrac12,\tfrac52,\tfrac25,{\tfrac {13}{5}},{\tfrac {5}{13}},{\tfrac {34}{13}},
{\tfrac {13}{34}},{\tfrac {89}{34}},{\tfrac {34}{89}},{\tfrac {233}{89}},{
\tfrac {89}{233}},{\tfrac {610}{233}},{\tfrac {233}{610}},{\tfrac {1597}{
610}},\dots)\\
&=
S(x; d_0, d_1, d_2,\dots).
\end{align*}
\end{thm}

To prove Theorem~\ref{t-1}, we define $R_n$ for $n\ge-1$ by
\begin{align*}
R_{-1} &= \frac32 - \frac12\sqrt{\frac{1-5x}{1-x}},\\
R_{2n+1}&=d_{2n+1} + \frac{1-3x-\5}{2x},\quad n\ge 0,\\
R_{2n}&=\frac{d_{2n}}{1-x R_{2n+1}},\quad n\ge 0.
\end{align*}

One can easily check that $R_m$ is a power series in $x$ with constant term
$d_m$ (with $d_{-1}=1$), though this will follow from Lemma \ref{l-1}.

\begin{lem}
\label{l-1}
For $m\ge-1$, we have
\[R_m = \frac{d_m}{1-x R_{m+1}},\]
where $d_{-1}=1$. 
\end{lem}
\begin{proof}
By definition, this is true if $m$ is even. 
Thus it is enough to prove that for $n\geq0$, 
\[R_{2n-1} = \frac{d_{2n-1}}{1-x R_{2n}} = \frac{d_{2n-1}}{1-\displaystyle\frac{x d_{2n}}{1-x R_{2n+1}}},\]
which is equivalent to
\begin{equation}
  \label{eq:2}
R_{2n+1} = \frac{1}{x} - \frac{1}{d_{2n-1}} - \frac{1}{R_{2n-1}-d_{2n-1}}.  
\end{equation}
We can check \eqref{eq:2} directly for $n=0$. Assume $n\geq1$. 
Then the right-hand side of \eqref{eq:2} is equal to
\begin{align*}
\frac{1}{x} - \frac{1}{d_{2n-1}} - \frac{2x}{1-3x-\5}
&=  \frac{1}{x} - \frac{1}{d_{2n-1}} - \frac{2x \left(1-3x +
  \5 \right)}{(1-3x)^2-(1-6x+5x^2)}\\
&=  \frac{1}{x} - \frac{1}{d_{2n-1}} - \frac{1-3x +
  \5}{2x}\\
&=  3- \frac{1}{d_{2n-1}} + \frac{1-3x -
  \5}{2x}.
\end{align*}
Since $3- 1/d_{2n-1} = d_{2n+1}$ by Lemma~\ref{thm:2}, we are done.
\end{proof}

\begin{proof}[Proof of Theorem \ref{t-1}]
It follows from Lemma \ref{l-1} that 
\begin{align*}
\frac32 - \frac12\sqrt{\frac{1-5x}{1-x}}&=R_{-1}=\frac{1}{1-xR_0}=\cfrac{1}{1-
  \cfrac{d_0x}{1-xR_1}}   =\cfrac{1}{1-
  \cfrac{d_0x}{1-
   \cfrac{d_1x}{1-xR_2}}} =\cdots.
\end{align*}
Continuing, and taking a limit, gives the S-fraction for $R_{-1}$.
\end{proof}

By \eqref{eq:gf}, Theorem~\ref{t-1} and Corollary~\ref{thm:4}, we obtain the
following which proves Theorem~\ref{thm:main}.
\[
\sum_{n\geq0} \#\NC_2(n) x^n = \frac32 -\frac12 \sqrt{\frac{1-5x}{1-x}}
=\sum_{n\geq0} \dyck_n(d) x^n 
=\sum_{n\geq0} \mot_n(b,\lambda) x^n
\]

\section{A combinatorial proof}

Let $b$, $\lambda$ and $d$ be the sequences defined in \eqref{eq:3} and
\eqref{eq:7}.

Recall that in the previous section we have shown that
$\dyck_n(d)=\mot_n(b,\lambda)$ by changing a Dyck path of length $2n$ to a
Motzkin path of length $n$. We can do the same thing after deleting the first
and the last steps of a Dyck path. More precisely, for a Dyck path of length
$2n$, we delete the first and the last steps, take two steps at a time in the
remaining $2n-2$ steps, and change $UU$, $UD$, $DU$ and $DD$, respectively, to
$U$, $H$, $H$ and $D$. Then we obtain a Motzkin path of length $n-1$.
This argument shows that
\begin{equation}
  \label{eq:1}
\dyck_n(d)=d_0\cdot \mot_{n-1}(\aa,\bb)=\mot_{n-1}(\aa,\bb),
\end{equation}
where $\aa_n = d_{2n}+d_{2n+1}$ and $\bb_n = d_{2n+1}d_{2n+2}$. By \eqref{eq:7}
and Lemma~\ref{thm:2}, we have $\aa=(2,3,3,\ldots)$ and $\bb=(1,1,\ldots)$. 
Note that we can also prove Theorem~\ref{t-1} using \eqref{eq:1}.

To find a connection between $\mot_{n-1}(\aa,\bb)$ and $\NC_2(n)$ we need the
following definition. 

A \emph{\sch. path} of length $2n$ is a lattice path from $(0,0)$ to $(2n,0)$
consisting of up steps $U=(1,1)$, down steps $D=(1,-1)$ and double horizontal
steps $H^2=HH=(2,0)$ that never goes below the $x$-axis.  Let $\scheven(n)$
denote the set of \sch. paths of length $2n$ such that all horizontal steps have
even height.

\begin{prop}\label{thm:5}
  Let $\aa=(2,3,3,\ldots)$ and $\bb=(1,1,\ldots)$. Then, for $n\geq1$, we have
\[\mot_{n}(\aa,\bb) = \#\scheven(n).\]
\end{prop}
\begin{proof}
  From a Motzkin path of length $n$ we obtain a \sch. path of length $2n$ as
  follows. Change $U$ and $D$ to $UU$ and $DD$ respectively. For a horizontal
  step $H$, if its height is $0$, we change it to either $UD$ or $HH$, and if
  its height is greater than $0$, we change it to either $UD$, $DU$ or
  $HH$. Then we get an element of $\scheven(n)$. Since the weight of a
  horizontal step $H$ in the Motzkin path is equal to the number of choices, the
  theorem follows.
\end{proof}

\begin{remark}
  The definition of $\scheven(n)$ in \cite{Kim2009b} is the set of \sch. paths
  of length $2n$ which have no peaks at even height. From such a path, by
  changing all the horizontal steps at odd height to peaks, we get a \sch. path
  whose horizontal steps are all at even height, and this transformation is
  easily seen to be a bijection.
\end{remark}

Kim \cite{Kim2009b} found a bijection between $\NC_2(n)$ and $\scheven(n-1)$.
Using Kim's bijection in \cite{Kim2009b}, Proposition~\ref{thm:5}, \eqref{eq:1} and
Corollary~\ref{thm:4} we finally get the following sequence of identities which
implies Theorem~\ref{thm:main}:
\[\#\NC_2(n) = \#\scheven(n-1) = \mot_{n-1}(\aa,\bb)
=\dyck_n(d)=\mot_n(b,\lambda). \]

\section*{Acknowledgement}
The authors would like to thank the anonymous referees for their careful reading
and helpful comments which significantly improved the presentation of the paper.

\bibliographystyle{plain}

\end{document}